\documentclass[a4paper]{article} 
\usepackage{amsmath,amsfonts,amssymb,amscd}
\usepackage[all]{xy}

\newtheorem{thm}{Theorem}[section]

\newtheorem{rema}[thm]{Remark}

\newtheorem{exx}[thm]{Example}

\def\A{\mathcal A}
\def\del{\partial}

\def\goesto{\mapsto}
\def\dcov{d^\Na}

\def\gr{gr}
 
\def\delstar{\pa^*}

\def\ST{\mathcal{S}}
\def\GL{{\rm GL}}
\def\SL{{\rm SL}}
\def\Ups{\Upsilon}

\def\rr{{\mathbb{R}}}

\let\i=\iota
\let\Na=\nabla

\def\embed{\hookrightarrow}
\def\bdot{\bullet}

\def\fsl{{\mathfrak{sl}}}
\def\gl{{\mathfrak{gl}}}

\newcommand{\newc}{\newcommand}
\newcommand{\cov}{\nabla}

\newcommand{\id}{\operatorname{id}}
\newcommand{\Ker}{\operatorname{Ker}}
\renewcommand{\Im}{\operatorname{Im}}
\newcommand{\Imm}{\operatorname{Im}}
\newcommand{\Idd}{\operatorname{Id}}

\let\ccdot\cdot
\def\cdot{\hbox to 2.5pt{\hss$\ccdot$\hss}}

\newcommand{\be}{\beta}
\newcommand{\ga}{\gamma}
\newcommand{\de}{\delta}

\newcommand{\ka}{\kappa}

\newcommand{\om}{\omega}
\renewcommand{\phi}{\varphi}
\newcommand{\ph}{\varphi}

\newcommand{\si}{\sigma}

\newcommand{\La}{\Lambda}
\newcommand{\Ga}{\Gamma}

\newcommand{\Om}{\Omega}

\newc{\aI}{\mbox{\boldmath{$ I$}}}
\newc{\aR}{\mbox{\boldmath{$ R$}}}
\newc{\aDeR}{\mbox{\boldmath{$ U$}}_B{}^P{}_C{}^Q}
\newc{\al}{\mbox{\boldmath$ \Delta$}}             
\newc{\nda}{\mbox{\boldmath$ \nabla$}}
\newc{\ad}{\mbox{\boldmath$ d$}}
\newc{\da}{\mbox{\boldmath$ \delta$}}
\newc{\aK}{\mbox{\boldmath{$ K$}}}
\newc{\aL}{\mbox{\boldmath{$ L$}}}

\newtheorem{theorem}{Theorem}[section]
\newtheorem{lemma}[theorem]{Lemma}

\newtheorem{corollary}[theorem]{Corollary}

\newcommand{\g}{{\mathfrak g}}

\newcommand{\pj}{{\bar{\i}}}
\newcommand{\pk}{{\bar{k}}}
\newcommand{\cA}{{\mathcal A}}

\newcommand{\cC}{{\mathcal C}}
\newcommand{\calC}{{\mathcal C}}

\newcommand{\cg}{{\mathcal G}}

\newcommand{\ce}{{\mathcal E}}
\newcommand{\cE}{{\mathcal E}}

\newcommand{\bW}{{\mathbb W}}

\newcommand{\bV}{{\mathbb V}}

\newcommand{\one}{{\mathbb I}}

\newcommand{\lapl}{\Box}

\newcommand{\Rho}{{\mbox{\sf P}}}

\newcommand{\End}{\operatorname{End}}
\newcommand{\im}{\operatorname{im}}

\newcommand{\Id}{\operatorname{Id}}

\newcommand{\bN}{{\Bbb N}}

\let\i=\iota

\let\G=\Gamma

\newc{\strutdd}{\rule{0mm}{5mm}}

\usepackage{ifthen}

\newc{\tensor}[1]{#1}

\newc{\Mvariable}[1]{\mbox{#1}}

\newc{\down}[1]{{}_{
\ifthenelse{\equal{#1}{;}}{|}{#1}}}

\newc{\up}[1]{{}^{#1}}
\newc{\C}{C}

\newc{\JulyStrut}{\rule{0mm}{6mm}}
\newc{\midtenPan}{\mbox{\sf S}}
\newc{\midten}{\mbox{\sf T}}
\newc{\midtenEi}{\mbox{\sf U}}
\newc{\ATen}{\mbox{\sf E}}
\newc{\BTen}{\mbox{\sf F}}
\newc{\CTen}{\mbox{\sf G}}

\def\sideremark#1{\ifvmode\leavevmode\fi\vadjust{\vbox to0pt{\vss
 \hbox to 0pt{\hskip\hsize\hskip1em
 \vbox{\hsize3cm\tiny\raggedright\pretolerance10000
 \noindent #1\hfill}\hss}\vbox to8pt{\vfil}\vss}}}%

\input epsf
\input xy
\xyoption{all}

\newtheorem{definition}[theorem]{Definition}

\newtheorem{remark}[theorem]{Remark}

\newcommand{\bd}{\begin{definition}}
\newcommand{\ed}{\end{definition}}

\newcommand{\br}{\begin{remark}}
\newcommand{\er}{\end{remark}}
\newcommand{\gog}{{\mathfrak g}}
\newcommand{\gop}{{\mathfrak p}}

\newcommand{\bt}{\begin{tabular}}
\newcommand{\et}{\end{tabular}}

\renewcommand{\ker}{\operatorname{ker}}
\newcommand{\Hom}{\operatorname{Hom}}
\renewcommand{\Im}{\operatorname{Im}}

\def\cdot{\hbox to 2.5pt{\hss$\ccdot$\hss}}

\def\gog{\mathfrak{g}}

\def\gop{\mathfrak{p}}

\newcommand{\pa}{\partial}

\def\C{\mathbb{C}}

\def\cC{\mathcal{C}}

\def\cE{\mathcal{E}}
\def\cG{\mathcal{G}}

\def\al{\alpha}
\def\be{\beta}
\def\ga{\gamma}
\def\de{\delta}

\def\et{\eta}

\def\ka{\kappa}

\def\rh{\rho}
\def\si{\sigma}
\def\ta{\tau}
\def\ph{\varphi}

\def\om{\omega}
\def\Ga{\Gamma}

\def\La{\Lambda}

\def\Om{\Omega}

\def\na{\nabla}
\def\pa{\partial}

\def\ra{\rightarrow}

\def\sideremark#1{\ifvmode\leavevmode\fi\vadjust{\vbox to0pt{\vss
 \hbox to 0pt{\hskip\hsize\hskip1em
 \vbox{\hsize3cm\tiny\raggedright\pretolerance10000
  \noindent #1\hfill}\hss}\vbox to8pt{\vfil}\vss}}}%

\begin{document}

\title{On a new normalization for tractor covariant derivatives}


\author{Matthias Hammerl, Petr Somberg, Vladimir Sou\v cek, Josef \v Silhan
\\ \\ \\
\small M. H.: Fakult\"{a}t fuer Mathematik, Nordbergstra{\ss}e 15, 1090 Wien, Austria,\\ 
\small              matthias.hammerl@univie.ac.at\\           
\small P. S., V. S.: Mathematical Institute,
Sokolovska 83, Karlin, 180 00, Prague 8, Czech Republic, \\
\small somberg@karlin.mff.cuni.cz, {soucek@karlin.mff.cuni.cz}\\
\small J. \v S.: Masaryk University, Janackovo nam. 2a, 662 95 Brno, Czech Republic, \\
\small              {silhan@math.muni.cz}\\
}

\maketitle

\pagestyle{myheadings}

\begin{abstract}
A regular normal parabolic geometry of type $G/P$ on a manifold $M$
gives rise to sequences $D_i$ of invariant differential operators, known as
the curved version of the BGG resolution. These sequences are
constructed from the normal covariant derivative $\na^\om$ on the
corresponding tractor bundle $V,$ where $\om$ is the normal Cartan connection.
The first operator $D_0$ in the sequence is overdetermined
and it is well known that $\na^\om$ yields the prolongation of this operator
in the homogeneous case $M = G/P$.
Our first main result is the curved version of such a prolongation. This requires
a new normalization $\tilde{\na}$ of the tractor covariant derivative on $V$.
Moreover, we obtain an analogue for higher operators $D_i$. In that case one needs
to modify the exterior covariant derivative $d^{\na^\om}$ by differential
terms. Finally we demonstrate these results on simple examples in projective
and Grassmannian geometry. Our approach is based on standard techniques
of the BGG machinery.
\end{abstract}



\section{Introduction}
Let $G$ be a (real) semisimple Lie group and $P$ its parabolic subgroup. Following ideas of \'E. Cartan, the homogeneous space
$G/P$ is a flat model for a curved parabolic geometry of type $(G,P),$ which is specified
by  a couple
$(\cG,\om),$  where $\cG\ra M$  is a principal $P$-bundle and $\om$ is a Cartan connection. 
It is well known that such a geometry can be characterized by an underlying geometric   
structure on the manifold $M,$ together with
a suitable normalization condition (for more details, see \cite{cap-slovak-par}). 
Distinguished examples of  this procedure are the normal Cartan connections
constructed for a conformal structure by \'E. Cartan and for a CR structure by Chern 
and Moser.
Let us consider a regular normal parabolic geometry $(\cG,\om)$ of type $(G,P).$
For any
$G$-module $\bV,$ the tractor bundle $V$ over $M$ is (by definition) the vector bundle 
associated to $\cG$ and the representation $\bV$ (restricted to $P$).
The normal Cartan connection $\om$ on $\cG$ then induces the tractor
covariant derivative $\Na^{\om}$
on $V,$ which  is then used in various problems in analysis
and/or geometry on $M$ (e.g., for constructions of differential 
invariants on the corresponding parabolic geometry). For example, it plays the key role
in the construction
of the Bernstein-Gel'fand-Gel'fand (BGG) sequences of invariant differential operators (see \cite{CSS,CD}) 
and prolongation procedures for first operators in the BGG sequences (see, e.g. \cite{BCEG}).
 
In particular, there is a lot of interest in the study of properties of the first operators in the BGG sequences, or their semilinear version. Ideas behind the construction of these operators by the BGG machinery
can be helpful in such problems. The construction uses  tractor covariant derivatives
acting on tractor bundles  and suitable splitting operators
(for details, see Sect. 3). In some simple cases, there is a one-to-one correspondence between solutions
of the first BGG equation and the kernel of the corresponding tractor covariant derivative. In other words, the
tractor covariant derivative is  the prolongation of the first BGG operator.  But such a
simple correspondence between solutions of the first BGG equation and the kernel of the tractor
covariant derivative is not valid in general. 

A general scheme for a prolongation of the first BGG operator (and its semilinear version)
was introduced in \cite{BCEG}, for a generalization to the contact cases, see \cite{eg}. The procedure used 
in \cite{BCEG} is efficient but not invariant.
In quite a few special cases (see \cite{Cap,GS,eastwood-matveev,E, DT}), several authors found an invariant way how to
compute a deformation of the normal tractor covariant derivative having the property that its
kernel can be identified with solutions of the first BGG sequence. 

The new normalization of 
tractor covariant derivatives developed in the paper is motivated by a wish to extend these examples
to a general scheme.
We shall study the problem of a suitable normalization for tractor covariant derivatives
for a general parabolic geometry in a systematic way and show that there is  a
distinguished alternative of the usual normalization
of tractor covariant derivatives on tractor bundles giving directly a canonical prolongation 
of the first BGG operator in an invariant way.

The normal  tractor covariant derivative is induced from the normal
Cartan connection on
the principal bundle $\cg.$ An important observation is that if we want to
find a covariant derivative on tractor bundles giving the invariant prolongation of the 
first BGG operator, it is necessary to adapt (in contrast to $\nabla^\omega$) the normalization condition 
to a choice of the tractor bundle under consideration.

The main results of the paper can be described as follows. Let us consider  
a regular normal parabolic geometry of type $(G,P)$ given by the couple $(\cG,\om).$ 
For any irreducible $G$-module $\bV,$  there is the associated
covariant derivative $\na^\om$   on the  associated vector bundle $V.$ The space of all
covariant derivatives on $V$ is the affine space modelled on the vector space
$\cE^1(\End V).$ We want to find a deformation of $\na^\om$  by $\Phi\in\cE^1(\End V)$ 
 satisfying a new normalization condition (adapted to the choice of $\bV$) in such a way that
the resulting covariant derivative will have suitable properties.

The deformation $\Phi$ cannot be chosen arbitrarily. Firstly, the construction of the BGG sequence leads to the requirement to preserve
the lowest homogeneous component of $\na^\om$ (having homogeneity zero), hence we
shall restrict to $\Phi\in(\cE^1(\End V))^1,$ where the superscript $1$ indicates
that $\Phi$ should have the (total) homogeneity bigger or equal to one.
The aim to have good properties of the new covariant derivative in
the prolongation procedure for the first BGG operator induces further restrictions
on a choice of $\Phi.$ They will be expressed by properties of values of $\Phi(s)\in
\cE^1(V),$ where $s$ is a section of $V.$
It leads to the following class of covariant derivatives on the tractor bundle $V.$
 
\begin{definition}
Let $\om$ be the regular normal Cartan connection on the principle bundle $\cg$ and let $\na^\om$ be
the associated covariant derivative on the associated vector bundle $V.$
The class $\cC$ of admissible covariant derivatives on $V$ is defined by
$$
\cC=
\left\{
\na=\na^\om+\Phi|\Phi\in \Im (\pa^*_V\otimes \Id_{V^*}),\Phi\in(\cE^1(\End V))^1,  
\right\},
$$
where $\pa^*_V$ is the Kostant differential  corresponding to homology
 of  $\g_-$ with values in $\bV.$
 
The condition $\Phi\in \Im (\pa^*_V\otimes \Id_{V^*})$ is equivalent to the property
$\Phi(s)\in  \Im\;\pa^*_V\subset \cE^1(V)$ for all $s\in\Gamma(V),$  
where $\Gamma (V)$ denotes the space of sections of $V.$
\end{definition}

The main theorem of the paper is then

\begin{theorem}
There exists a unique covariant derivative $\na\in\cC$ with the property
$$
(\pa^*_V\otimes \Id_{V^*})(R^\na)=0,
$$
where $R^\na\in\cE^2(\End V)$ is the curvature of $\na.$

Again,   the condition $(\pa^*_V\otimes \Id_{V^*})(R^\na)=0$ can be equivalently
  expressed as
  the condition $\pa^*_V(R^\na(s))=0$ for all sections $s$ of $V.$

\end{theorem}

The new covariant derivative $\Na$ constructed in Theorem 1.2  gives a prolongation of the first
BGG operator, hence we shall call the covariant derivative satisfying this new normalization condition
{\it the prolongation covariant derivative.} 
 The next main result is the theorem stating this property.
 
\begin{theorem}
Let us consider a parabolic geometry $(\cG,\om)$ modeled on a couple $(G,P).$
There is a one-to one correspondence between the kernel of 
the first BGG operator for a $G$-module $\bV$ and the kernel of the
prolongation covariant derivative on the associated bundle $V$ over $M.$
\end{theorem}

In the second part of the paper, we extend the previous construction to other operators in the BGG sequence.
In these cases, we have to consider a more general deformation of the exterior derivative
$d^\Na$ by adding a differential term (instead of just an algebraic one, which
was sufficient for the first operator in the BGG sequence).

Finally, we compare the general procedure developed in the paper with particular results obtained
in some special cases and compute some other examples of the prolongation covariant derivatives.
They  come from projective and Grassmann geometry.

\section{Normalization of tractor covariant derivatives}
  
\subsection{The double filtration on $\End \bV$}
Let $G$ be a semisimple Lie group (real or complex) and $P$ its parabolic subgroup.
The choice of $P$ induces the grading  $\gog=\oplus_{i=-k}^k\gog_i$ 
on the Lie algebra of $G.$ 
Let $\bV$ be an irreducible module for $G.$
There is the grading element $E$ in $\gog_0$ acting by $i$ on $\gog_i.$ It
 can be additively shifted to $E'$ in such a way that eigenvalues of $E'$
on $\bV$  are integers between $0$ and $r$ for a suitable positive integer $r.$
Eigenvalues of $E'$ on $\bV^*$ are then integers between $-r$ and $0.$ Then we get decompositions of $\bV,$ resp. $\bV^*,$
into the corresponding eigenspaces 
$$\bV=\oplus_{\pj=0}^r\bV_\pj,\;\;\;\;\ 
\bV^*=\oplus_{\pj'=-r}^{0}\bV_{\pj'}^*.
$$  A similar decomposition of $\gog_+$ is given by
$\gog_+=\gog_1\oplus\ldots\oplus\gog_k.$

A bigrading on $\End \bV\simeq \bV\otimes \bV^*$ is then given by
 $$
\End\bV=
\oplus_{\pj=0}^r\oplus_{\pj'=-r}^0(\bV_{\pj}\otimes\bV^*_{\pj'}).
$$ 
Consequently there is a 'diagonal' grading on $\End\bV$ given by
$$
\End\bV=\oplus_{\ell=-r}^r (\End\bV)_{\ell};\;\;\;
(\End\bV)_{\ell}:=  \oplus_{\pj+\pj'=\ell}\bV_{\pj}\otimes\bV^*_{\pj'},
$$
and 'vertical' (resp. 'horizontal') gradings given by
$$
\End\bV=\oplus_{\pj=0}^r (\End\bV)_{\pj};\;\;\;
(\End\bV)_{\pj}:=\bV_{\pj}\otimes\bV^*
$$
$$
\End\bV=\oplus_{\pj'=-r}^0(\End\bV)_{\pj'};\;\;\;
(\End\bV)_{\pj'}:=\bV\otimes\bV^*_{\pj'}.
$$
The diagonal grading is independent of the normalization of the grading of $\bV.$
In what follows, we shall use the diagonal and the vertical gradings on $\End\bV.$

The gradings are not $P$-invariant. 
We shall hence consider filtrations induced by gradings above.
For the diagonal grading, we shall define the filtration  by a choice of subspaces
$$
(\End\bV)^{\ell}=\oplus_{k\geq\ell} (\End\bV)_{k},
$$
while for the horizontal grading, the filtration is defined by
$$
(\End\bV)^\pj=\oplus_{\pk\geq\pj} (\End\bV)_{\pk}.
$$
The grading of $\gog_+$ also gives the standard filtration 
$\gog^k\subset\ldots\subset\gog^1=\gog_+.$

These filtrations (together with the filtration on $\gog_+$) induce also the filtrations on the chain spaces
$\Lambda^j(\gog_+)\otimes \End\bV$ for the 
Lie algebra homology  and cohomology complexes.
The differentials  in the Lie algebra
(co)ho\-mo\-logy of $\gog_-$ with values in $\gog$-modules $\bW$ 
are the maps 
$\pa_{\bW}:\Lambda^j(\gog_+)\otimes \bW\mapsto \Lambda^{j+1}(\gog_+)\otimes \bW\,$
 resp. 
 $\pa^*_{\bW}:\Lambda^j(\gog_+)\otimes \bW\mapsto \Lambda^{j-1}(\gog_+)\otimes \bW.$ 
If $\bW=\End \bV\simeq \bV\otimes\bV^*$ for a $\gog$-module $\bV,$ we shall 
denote operators $\pa_\bV\otimes \Id_{\bV^*},$ resp. $\pa^*_{\bV}\otimes \Id_{\bV^*},$ 
simply by  $\pa_\bV,$ resp. $\pa_{\bV}^*.$ It should not lead to any confusion.

The definition of operators $\pa_\bV$ and $\pa_\bV^*$ implies immediately
that they preserve both horizontal and diagonal gradings on $\Lambda^j(\gog_+)\otimes\End \bV.$ 
Hence they respect
both horizontal and diagonal filtrations  on  $\Lambda^j(\gog_+)\otimes\End \bV.$
We shall use bellow the  induced operators  between the graded bundles associated to the horizontal filtration and we shall denote them by $gr\,\pa_{\bV},$ resp. $gr\,\pa_{\bV}^*.$

\subsection{Induced operators on associated graded bundles}
The spaces of $j$-forms on $M$ with values in a bundle $W$  will be
denoted by $\cE^j(W).$ They are isomorphic to the bundle induced  by the $P$-module $\Lambda^j(\gog_+)\otimes\bW.$ Similarly, the tangent bundle is isomorphic to the bundle
associated to the $P$-modul $\gog/\gop.$
All filtrations mentioned above are $P$-invariant and they consequently induce the corresponding
filtrations on $\cE^j(\End V).$
We shall need, in particular, the diagonal filtrations 
$(\cE^j(\End V))^\ell,$
resp. the vertical filtration $(\cE^j(\End V))^\pj,$ induced on $\cE^j(\End V).$
We shall denote by $gr_\ell(\cE^j(\End V)),$ resp. $gr_\pj(\cE^j(\End V))$ the associated
graded bundles.

The operator $gr\pa_\bV^*$ and $gr\pa_\bV$ are $P$-equivariant, hence they induce  well defined maps $\pa_V^*, $ resp. $\pa_V,$ between
the corresponding associated graded bundles. 
 
We shall denote by $gr\,\pa_V,$ resp. $gr\,\pa_V^*,$ the direct sum of all maps
$gr_\pj\,\pa_V,$ resp. $gr_\pj\,\pa_V^*$ acting on the direct sum
$gr \, \cE^j(\End V):= \oplus_\pj gr_\pj(\cE^j(\End V)).$
The operators $gr\,\pa_V,$ and $gr\,\pa_V^*$ have then usual properties of the Kostant differentials.
In particular, they are dual to each other (with respect to a suitable scalar product), which implies
usual properties of their kernel and images (the Hodge decomposition).

 Note also that $\cE^j(V)\otimes V^*=\cE^j(\End V ).$ Hence the standard filtration
on $\cE^j(V)$ is transferred (by the tensor product with $V^*$) to the
horizontal grading on $\cE^j(\End V).$
As an immediate corollary, we get that $\phi\in\cE^j(\End V)^\pj$ if and
only if $\phi s\in  \cE^j(V)^\pj$ for all sections $s\in\cE^0(V).$

\subsection{A choice of normalization}
Let us consider a  regular parabolic geometry $(\cg,\om)$ over $M$ with the homogeneous model
given by a couple $(G,P).$ For an irreducible $G$-module $\bV,$ we shall consider the associated tractor bundle $V$ on $M.$ The curvature $\kappa$ of the  Cartan connection $\om$  is a two-form with values in the adjoint tractor bundle ${\mathcal A}\simeq \mathcal G\times_P\gog.$
The usual normalization condition for $\omega,$
expressed in terms of the Kostant differential $\pa^*$ corresponding to homology
 of  $\g_-$ with values in $\g,$  requires the curvature $\kappa$ to be $\pa^*$-closed. 
In terms of an associated covariant derivative $\Na^\om$ on $V,$ the curvature $R^{\Na^\om}$ of $\Na^\om$
is a two-form with values in $\End V$ and the normalization condition can be expressed using
the Kostant differential $\pa^* $ for $\End V $  as 
$$
\pa^*(R^{\Na^\om})=0.
$$

 Given a choice of the bundle $V,$ we are going to change the normalization
condition for a covariant derivative $\Na$ on $V.$
  Let $\Idd_{V^*}$ denote the identity map on $V^*.$
As above in the algebraic version, we shall consider operators
$$
\pa_V\otimes \Idd_{V^*},\; \pa_V^*\otimes \Idd_{V^*}.
$$
acting on forms $\cE^j(\End V)$ with values in $\End V\simeq V\otimes V^*.$
Abusing the notation, we shall denote them by $\pa_V,$ resp. $\pa_V^*.$
It will always be clear whether the differentials act on forms with values
in $V$ or forms with values in $\End V.$

We shall now introduce a new normalization for covariant derivatives on $V.$

\begin{definition}
We shall call a covariant derivative $\Na\in\cC$ {\rm the prolongation covariant derivative},
if
$$
\pa_V^*(R^\Na)=0,
$$ 
where  $R^\Na\in\cE^2(\End V)$ is the curvature of $\Na.$
\end{definition}

The choice of its name should suggest that the new normalization condition gives better
properties to $\Na$ in the prolongation procedure for the first operator in the BGG sequence
corresponding to the representation $\bV$ (or its semilinear versions).

We shall need the following property.
\begin{lemma}\label{product}
If $\phi\in(\cE^1(\End V))^\pj$ and $\tau\in\cE^1(V),$ then
$$
\phi\wedge\tau\in(\cE^2(V))^{\pj+1}.
$$
\end{lemma}

\noindent
{\bf Proof.}

Indeed, we can decompose $\phi$ into homogeneous components
$$
\phi=\sum_j\al_j\otimes v_j\otimes w_j,\al_j\in\cE^1,v_j\in V,w_j\in V^*,
$$
where the sum of homogeneities of $\al_j$ and $v_j$ is greater or equal to $\pj.$ If we also decompose
$\tau$ as 
$$
\tau=\sum_k\beta_k\otimes u_k,\beta_k\in\cE^1,u_k\in V,
$$
then the expression
$$
\phi\wedge\tau=\sum_{j,k}w_j(u_k)\al_j\wedge\beta_k\otimes v_j
$$
clearly has summands of homogeneity greater or equal to $\pj+1.$
 \hfill$\square$

\subsection{The main lemma.}  

The key information for the normalization procedure is the following fact concerning
the induced change of the curvature.
\begin{lemma}\label{change}
Let $\Na_1,$ resp. $\Na_2,$ be two  covariant derivatives from $\cC$ related to each other
by 
the deformation 
$\Phi=\Na_2-\Na_1\in(\cE^1(\End V ))^1$
and let $R_1,$ resp. $R_2,$ be the corresponding curvatures.

If $\Phi\in(\cE^1(\End V ))^\pj,$ then $R_2-R_1\in(\cE^2(\End V ))^\pj$
and
$$
gr_\pj(R_2-R_1)=(gr_\pj\,\pa_V)(gr_\pj\,\Phi).
$$
\end{lemma}

\noindent
{\bf Proof.}

Let $\om$ be the normal Cartan connection for the chosen parabolic geometry and $\Na$ its associated
covariant derivative. It is well known that $\Na$ and $d^{\Na}$ preserve the standard
filtration on $\cE^j(V)$ and that the corresponding graded version of $\Na,$ resp.
$d^{\Na}$ is equal to $gr\,\pa_V.$ A shift of $\Na$ by $\Phi\in(\cE^1(\End V ))^1$
does not change this property, the same being true for $d^{\Na+\Phi}.$

 The change in the curvature is then
$$R_2-R_1=d^\Na\Phi+[\Phi,\Phi].
$$
The result clearly belongs to $\cE^2(\End V )^\pj,$ because the operator $d^\Na$
preserves the filtrations and we can use Lemma \ref{product} for the second term.

Then we get for  any $s\in\cE^0(V),$
$$
gr_\pj((d^\Na\Phi+[\Phi,\Phi])s)=
gr_\pj((d^\Na\Phi) s)=
$$
$$=gr_\pj(d^\Na(\Phi s)-\Phi\wedge (\Na s))=
gr_\pj(\pa_V(\Phi(s)))=
$$
$$
=(gr_\pj\pa_V)( gr_\pj(\Phi(s))).
$$ 
 \hfill$\square$
 
\subsection{Existence}

\begin{lemma}
Suppose that there is a  tractor covariant derivative $\Na\in\cC$ with the property
$$
\pa_V^*(R^\Na)\in\cE^1(\End V)^\pj,
$$
where $\pj$ is a number between $0$ and $r.$

Then there exists $\Phi\in\cE^1(\End V)^{1}\cap\cE^1(\End V)^\pj$ such that for 
$\tilde{\Na}=\Na+\Phi,$  we have
$$
\pa_V^*(R^{\tilde{\Na}})\in\cE^1(\End V)^{\pj+1}.
$$
\end{lemma}

\noindent{\bf Proof.}

The spaces $\{ \cE^1(\End V)^{1}\cap\cE^1(\End V)^\pj\}_{\pj=0}^r$ give a descending filtration of
the space $\cE^1(\End V)^1.$ The filtration is preserved by maps
$\pa_V$ and $\pa_V^*,$ hence they induce maps on the associated graded
bundle (we denote them for simplicity of notation by the same symbols as for the full
filtration of $\cE^1(\End V)$). 
The standard Kostant decomposition says that $\Ker gr\,\pa_V^*$ and $\Im gr\,\pa_V$ are complementary subspaces of the
graded bundle $gr\,\cE^1(\End V)^1.$ In particular, $gr\,\pa_V^*$ restricts to an isomorphism
of $\Im gr\,\pa_V$ to $\Im gr\,\pa_V^*.$

Hence there is an element $\phi\in gr_{\pj}(\cE^1(\End V)^1)$ such that
$$
(gr\,\pa_V^*)((gr\,\pa_V)(\phi))=gr_\pj(\pa_V^*(R^\Na)).
$$

We shall take any preimage $\Phi\in\cE^1(\End V)^{1}\cap\cE^1(\End V)^\pj$
of  $\phi$ and we shall define a corrected covariant derivative
by 
$\tilde{\Na}=\Na-\Phi.$

Due to Lemma \ref{change}, we get
\begin{eqnarray*}
gr_\pj(\pa_V^*(R^{\tilde{\Na}}))&=&
gr_\pj(\pa_V^*(R^{{\Na}}))-(gr\,\pa_V^*)(gr_\pj(R^\Na-R^{\tilde{\Na}}))=\\
&=& gr_\pj(\pa_V^*(R^{{\Na}}))-(gr\,\pa_V^*)((gr\,\pa_V)(gr_\pj (\Phi)))=0.
\end{eqnarray*}
Hence $\tilde{\Na}$ has the required properties.

\hfill $\square$

\begin{theorem}
For each irreducible $G$-module $\bV,$ 
there exists a prolongation covariant derivative $\Na\in\cC,$ i.e.,  we can find $\Na\in\cC$ such that
$$
\pa_V^*(R^\Na)=0.
$$
\end{theorem}

\noindent
{\bf Proof.}

The curvature function of the regular normal connection $\om$ for the corresponding parabolic geometry belongs
(by definition of regularity) to $\cE^2(\cA)^1,$ so 
$R^{\Na^\om}\in\cE^2(\End V)^1,$ and  $\pa_V^* ( R^{\Na^\om})
\in\cE^1(\End V)^1.$
Using Lemma 2.4., we get (by induction) the claim of the theorem.

\hfill $\square$

\subsection{Uniqueness.}

\begin{theorem}\label{uniq}
Suppose that $\Na_1,\Na_2$ are two covariant derivatives in $\cC,$  
both satisfying the normalization condition
$\pa_V^*(R^\Na)=0.$ Then 
$$
\Na_1=\Na_2.
$$
 \end{theorem}

\noindent
{\bf Proof.}

 Let $\Phi_1,\Phi_2\in\cE^1(\End V )^1\cap \Im \pa_V^*$ such that
$$
\Na_1=\Na^{\om}+\Phi_1;\;\Na_2=\Na^{\om}+\Phi_2.
$$
 
Denote by $R_1,$ resp. $R_2,$ the curvatures of $\Na_1,$ resp. $\Na_2.$
Then $\Phi=\Phi_2-\Phi_1$ belongs to $\cE^1(\End V )^1\cap \Im\pa_V^*.$
Suppose now that $\Phi\in\cE^1(\End V )^\pj.$
By assumption,
$gr_\pj(R_2-R_1)$ is in the kernel of $gr\, \pa_V^*.$
By Lemma \ref{change},  we have
$$
gr_\pj(R_2-R_1)=(gr\,\pa_V)(gr_\pj\,\Phi).
$$
But $\Ker gr \,\pa_V^*\cap \Im gr\, \pa_V$ is trivial hence
$gr_\pj(R_2-R_1)=0.$
Hence $gr_\pj\,\Phi$ is in the kernel of $gr\,\pa_V,$
and also in the image of $gr\,\pa_V^*,$ by assumption.
Hence $gr_\pj\,\Phi=0.$
By induction, $\Phi=0.$
\hfill$\square$

The construction above depends on some choices (e.g., a choice of a preimage $\Phi$
of $\phi$). Nevertheless, the uniqueness of the prolongation covariant derivative
shows that the result of the construction is independent of all choices. Hence we get
the following corollary.

\begin{corollary}
The prolongation covariant derivative is invariant. This means that it depends
only on the data of the chosen parabolic structure and the bundle $V.$
\end{corollary}

\section{The prolongation of the first BGG operator.}

The BGG complexes are sequences of invariant differential operators on a homogeneous model for
a given parabolic geometry.
A curved version of it, i.e., an extension of operators in the sequence to invariant differential operators
on general (non-flat) manifolds with a given parabolic structure 
 was first constructed in \cite{CSS}\ and the
construction was simplified and extended in \cite{CD}. The first operator
in such a sequence always gives an overdetermined system of invariant differential equations. A prolongation
of this operator was constructed for the case of $1$-graded parabolic geometries in \cite{BCEG}.
 However, the methods used there needed a choice of the Weyl structure, hence the resulting
covariant derivative was not invariant. We are now going to show that the normalization of tractor connections described in the paper
can be used to obtain invariant (natural) prolongations.

We begin by introducing the setting and basic operators of the BGG-machinery in
a generalized version needed for the next section.
Let $V$\ be a tractor bundle over $M$ with a covariant derivative $\na$ and the
exterior covariant derivative $\dcov:\ce^k(V)\goesto\ce^{k+1}(V).$ 
Recall from above that we have a well defined differential $\delstar=\pa^*_V:\ce^{k+1}(V)\to\ce^k(V)$.
The property $\delstar\circ\delstar=0$ allows us to define the cohomology $H_k$ as the
vector bundle quotient $H_k=\Ker\delstar/\Imm\delstar,$
   where $\Ker\delstar\subset\ce^k(V)$ is the space of cycles and $\Imm\delstar\subset\ce^k(V)$ is the space of boundaries.
The canonical surjection $\Ker\delstar_V\subset\ce^k(V)\goesto H_k$\ will be denoted by $\Pi_k$.


Due to regularity of the parabolic geometry under consideration, the operators $\dcov$\ are homogeneous of degree
zero with respect to the natural filtration of the spaces $\ce^k(V)$\ and they
induce the algebraic differential
$\gr\del_V:\gr(\ce^k(V))\goesto\gr(\ce^{k+1}(V))$\ on the associated
graded spaces. Thus it is possible to regard $\dcov$\ as a natural lift of $\gr(\del_V)$\ to a differential operator
from $\ce^k(V)$  to $\ce^{k+1}(V).$ 
 
The main ingredients in the BGG-machinery are the differential splitting
operators $L_k:H_k\goesto \Ker\delstar_V\subset\ce^k(V)$ with the property $\pa^*\circ d^\nabla \circ L_k=0.$  
This allows one to define the BGG-operators $D_k:H_k \goesto H_{k+1}$\ in 
the obvious way: $D_k:=\Pi_k\circ \dcov\circ L_k$. The definition is  encoded in the diagram  
\begin{align}\label{diagram1}
\xymatrix{
\ce^k(V)\ar[r]^{\dcov}&\ce^{k+1}(V)\\
\Ker\delstar\ar[u]^{i}  & \Ker\delstar\ar[u]^{i}\ar[d]^{\Pi_{k+1}}\\
H_k \ar[u]^{L_k}\ar[r]^{D_k} & H_{k+1} 
}
\end{align}
where $i$ denotes the inclusion.

We shall introduce the construction of the splitting operators in a more general situation, where the exterior covariant derivatives $d^\nabla$ on $\cE^k(V)$ will be substituted
by  general differential operators $E_k$ with suitable properties (see the theorem below).
The operators $D_k$ are defined by the same construction as the BGG operators and they
depend, in general, on the choice of $E_k.$ The theorem below shows that for  certain
classes of operators $E_k,$ the resulting operators $L_k$ and $D_k$ do not change.

\begin{thm}\label{BGG}
Let $(\ce^k(V))^j$ denote the filtration on $\ce^k(V)$ and let $gr(\ce^k(V))$ denote the associated
graded bundle, similarly for $\ce^{k+1}(V).$
Let $E_k$ be a  filtration preserving differential operator from  $\ce^k(V)$ to $\ce^{k+1}(V)$
with the property that the associated
graded map coincides with $gr\, \pa.$

Then for every $\si\in H_k,$ there exists a unique element $s\in\Ker\pa^*$ with  the following properties:

(1) $\Pi_k(s)=\si,$

(2) $E_k(s)\in\Ker\pa^*.$

Moreover, the mapping $L_k$ defined by  $\si\mapsto L_k(\si):=s$ is given by a differential operator.
The corresponding operator $D_k$ is then defined by
$$
D_k:=\Pi_{k+1}\circ E_k\circ L_k: H_k\mapsto H_{k+1}.
$$

Suppose that we change the operator $E_k$ to $\tilde{E}_k=E_k+\Phi_k,$ where the map $\Phi_k:\cE_k(V)\rightarrow\cE_{k+1}(V)$ is a differential
operator with values in $\Imm \pa^*$, and preserving the filtration
with the property that the associated
graded map is trivial.

Then the construction does not change the splitting operator $L_k$ and the  operator $D_k.$
\end{thm}

\noindent{\it Proof.}

The first part of the proof  follows the standard line of arguments.
The operator $\pa^*\circ E_k$ acts on $\cE^k(V)$ and it preserves $\Im \pa^*.$
It preserves the filtration and its graded version is, by assumption, given
by $gr(\pa^*)\circ gr(\pa),$ which is invertible on $\Im \pa^*.$
Hence also $\pa^*\circ E_k$ is invertible on $\Im \pa^*$ and it is possible
to show that its inverse $Q$ is a differential operator.

We can then define a differential operator $\hat{L}_k:=\Id-Q\circ\pa^*\circ E_k,$ which
restricts to zero on $\Im\pa^*.$  Hence it induces a well-defined differential
operator $L_k$ from $H_k$ to $\Ker\pa^*\subset\cE^k(V).$ It is easy to check that
the operator $L_k$ satisfies three properties 
$$
\Im L_k\subset \Ker \pa^*,\, \Pi_k\circ L_k=\Id,\, \pa^*\circ E_k\circ L_k=0.
$$

To show that $L_k$ is uniquely characterized by these properties, let us consider
$s_1,s_2\in \Ker\pa^*$ such that $E_k(s_i)\in\Ker\pa^*,\,i=1,2$ and $\Pi_k(s_1)=\Pi_k(s_2).$
 Then the difference $s=s_1-s_2$ belongs to $\Im\pa^*.$
 By definition of $\hat{L}_k,$ the relation $\pa^*\circ E_k(s)=0$ implies $\hat{L}_k(s)=s.$  On the other hand, $\hat{L}_k$ is trivial on $\Im\pa^*,$ Hence $\hat{L}_k(s)=0.$

To prove the last statement of the theorem, we shall consider  a  section $s$\ of $\ce^k(V)$. 
The new operator $\tilde{E}_k$ preserves the filtration and the induced graded map is still $\gr\,\pa.$ Since $(\tilde{E}_k -E_k)s$ belongs to
$ \Im\delstar_V$,
one has $\tilde{E}_k(s)\in \Ker\delstar_V$\ iff $E_k (s)\in \Ker\delstar_V$, which shows that $\tilde{L}_k=L_k$.
Thus, for $\si\in H_k$, one has $(\tilde{E_k}\, \tilde L_k -E_k\, L_k)\si\in \Im\delstar_V$,
but this lies in the kernel of the projection $\Pi_{k+1}:\Ker\pa^*\goesto H_{k+1}$. 
\hfill$\square$

Now we want to discuss the relation between $\Ker E_k$ and $\Ker D_k.$ For that, we have
to consider two consecutive operators $E_k$ and $E_{k+1}$ at the same time. They define
two splitting operators $L_k$ and $L_{k+1}.$ We get in such a way the diagram
\begin{align}\label{diagram2}
\xymatrix{
\cE^k(V)\ar[r]^{E_k} & \cE^{k+1}(V)\\
H_k \ar[u]^{L_k}\ar[r]^{D_k} & H_{k+1} \ar[u]^{L_{k+1}}
}
\end{align}
which, in general, does not commute but there is a convenient criterion for its
commutativity.

\begin{thm}
The diagram (\ref{diagram2}) commutes if and only if $\pa^*\circ E_{k+1}\circ E_k(s)=0$
for all sections $s\in\Im L_k\subset\ce^k(V).$
\end{thm}

\noindent {\it Proof.}

The values of $L_k$ are uniquely characterized by the conditions $L_k(\si)\in\Ker\pa^*$ and
$E_k\circ L_k(\si)\in\Ker\pa^*.$
Similarly, the values of $L_{k+1}$ are characterized by $L_{k+1}(\tau)\in\Ker\pa^*$ and
$E_{k+1}\circ L_{k+1}(\tau)\in\Ker\pa^*.$
Hence $E_k\circ L_k(\si)=L_{k+1}\circ D_k(\si)$ iff
$E_{k+1}\circ E_k\circ L_k(\si)\in\Ker\pa^*$ for all $\si\in H_k.$
 \hfill$\square$

If the diagram above is commutative, we get immediately a one-to-one correspondence between 
$\Ker E_k\cap\Ker\pa^*$and $\Ker D_k.$

\begin{thm}
Suppose that the diagram (\ref{diagram}) commutes.
Then $\Pi_k$\ and $L_k$\ restrict to inverse isomorphisms between $\Ker E_k\cap\Ker\pa^*$
 and  $\Ker D_k$.
\end{thm}

\noindent{\it Proof.}
  Let $s$\ be in $\Ker E_k\cap\Ker\pa^*.$ Then  $s=L_k(\Pi_k(s))$
  by definition of $L_k,$   and $\Pi_k(s)\in\Ker D_k$ by definition of
  $D_k.$
  
  On the other hand, if $D_k(\si)=0,$  then commutativity of the diagram implies
  that also
  $$
  L_{k+1}\circ D_k(\si)=E_k\circ L_k(\si)=0,
  $$ 
  hence 
  $L_k\in \Ker E_k\cap\Ker\pa^*.$
   
  And by definition of $L_k,$ we have $\Pi_k\circ L_k=\Idd.$
\hfill$\square$

\vskip 1mm
\noindent
Now we can return back to properties of the prolongation covariant derivative $\Na$ on $V.$
Using the above claims in the special case of the first square and operators
$E_0=\na$ and $E_1=d^\na,$ we see immediately that $E_1\circ E_0= R^\na.$
Hence we get the following corollary.

\begin{corollary}\label{cor-prolong}
Consider a tractor bundle $V$ and the corresponding prolongation covariant
derivative $\Na.$  Set $E_0=\Na $ and $E_1=d^\Na.$

Then the square constructed using these two operators commute and the covariant
 derivative $\na$ gives a prolongation of the first BGG operator $D_0.$
 In particular, the splitting operator $L_0$ induces a one-to-one correspondence
 between the space of parallel sections of $V$ with respect to $\na$ and the kernel
 of the first BGG operator $D_0.$
\end{corollary}

\vskip 2mm
\noindent
{\bf Remark.}
	In the case of a $1$-graded geometry, it was shown in \cite{BCEG}\ that
	the map $L_0:H_0\goesto V$\ induces an isomorphism of $J^k H_0$\ with $V_{\leq k}$\
	for every $k$\ such that the homology of $H_1(\g_-,\bV)$\ sits in
	homogeneity $>k$. Thus, for every operator $\tilde D_0:H_0\goesto H_1$\ which differs from the standard
  BGG-operator $D_0$\ by a linear differential operator of order $\leq k$,
  there is a map $\Psi\in\ce^1(\End V )$\ with values in $\Ker\delstar_V$\ such
  that its induced first BGG-operator coincides with $\tilde D_0$. 
  The mapping $\Psi$\ is unique up to maps $\ce^1(\End V )$\ with values in $\Im\delstar_V$, and
  it is thus easy to see that therefore the resulting normalized connection
  $\tilde\na=\na+\Psi+\Phi$\ doesn't depend on the choice of $\Psi$. Thus, natural deformations
  of $D_0$\ of low enough order can be prolonged naturally as well. We remark that
  a similar procedure works in the case of general graded parabolic geometries, where
  one has to use the filtration of the manifold for a suitable version of jet bundles.

\section{Prolongation covariant derivatives for the whole BGG sequence}

In this section we shall treat the problem  considered above in the case of other squares 
of the BGG sequence. We want to deform the exterior covariant derivative $d^\Na$
on $k$-forms in such a way that all squares in the generalized BGG construction
will commute, and, at the same time, the BGG operators $D_k$ will not change.
In fact, we shall succeed to keep both the BGG operators $D_k$ and the splitting operators  
$L_k$ unchanged. The deformation of $d^\Na$ on $\cE^k(V)$ will have, however, a different character.
It will be replaced by $E_k:=d^\Na+\Phi_k,$ where $\Phi_k$ is a linear differential operator
mapping $\cE^k(V)$ to $\cE^{k+1}(V).$  Hence the deformation $\Phi_k$ will not be, in general, algebraic.
Necessary tools were already prepared in the previous section (Theorems 3.1.- 3.3.). Methods described
in this section can also be applied to the first square but they give different answer (and also in this case the deformation $\Phi_0$ will not be algebraic in general).

To describe allowed deformations of the exterior derivative $d^\na,$ we shall introduce the following
notation.
There are two different filtrations on the space $A:=\Hom(\cE^k(V),\cE^{k+1}(V)).$
The diagonal filtration $A^j$ is induced by the standard filtration on $\cE^k(V),$ which is defined
by the condition $\Phi(s)\in(\cE^{k+1}(V))^{{a}+j}$ for all $s\in(\cE^k(V))^a.$
The other (vertical) filtration  $A^\pj$ is defined by the condition
$\Phi(s)\in(\cE^{k+1}(V))^\pj$ for all $s\in\cE^k(V).$
In this paragraph, we shall use symbols $\pa$ and $\pa^*$ for the Kostant differential
associated to the spaces $\cE^k(V).$
Recall that the class $\cC$ of admissible covariant derivatives on $V$ was defined by
$$
\cC=
\left\{
\na=\na^\om+\Phi|\Phi\in \Im (\pa^*_V\otimes \Id_{V^*}),\Phi\in(\cE^1(\End V))^1
\right\}.
$$

We shall consider the following spaces $\cC_k$ of deformations.

\begin{definition}
The space of allowed deformations will be defined by
$$
\cC_k:=\{E_k\in\Hom(\cE^k(V),\cE^{k+1}(V))| E_k=d^\na+\Phi,\,\Phi\in A^1, \Im \Phi\subset \Im \pa^*\}
$$
\end{definition}

\begin{theorem}  

{\ }
\noindent (1)
Let $\Na$ be any covariant derivative from $\calC.$ Let us consider the BGG sequence
with the splitting operators $L_k$ and the BGG operators $\{D_k\}$ induced (via Theorem \ref{BGG}) by operators $E_k=d^\na$ 
\begin{align}\label{diagram}
\xymatrix{
\cE^k(V)\ar[r]^{d^\nabla} & \cE^{k+1}(V)\\
H_k \ar[u]^{L_k}\ar[r]^{D_k} & H_{k+1} \ar[u]^{L_{k+1}}
}
\end{align}






Then there exists a collection of differential operators 
$\Phi_k\in\cC_k$ such that $\pa^*\circ d^\na\circ(d^\na+\Phi_k)=0.$
Moreover, the collection $\Phi_k$ with these properties is unique.

\vskip  2mm
\noindent
(2) As a consequence, the diagrams 
\begin{align}\label{BGG_k} 
\xymatrix{
\cE^k(V)\ar[r]^{d^\nabla+\Phi_k} & \cE^{k+1}(V)\\
H_k(V) \ar[u]^{L_k}\ar[r]^{D_k} & H_{k+1}(V) \ar[u]^{L_{k+1}}
}
\end{align}





 commute for all $k=0,1,\ldots,n-1$.
 
 If $\na$ depends only on data of the chosen parabolic geometry, the same is true
 for operators $E_k=d^\na+\Phi_k.$
 \end{theorem}  

\noindent{\it Proof.}
Let us choose $k=0,\ldots,n-1$
and  consider the square (\ref{BGG_k}) in  the generalized BGG sequence 
constructed using 
operators $d^\Na,$ where $\Na$ is any covariant derivative from $\calC.$   
We shall first prove the first assertion of the theorem.

The spaces $\{ A^{1}\cap A^\pj\}_{\pj=0}^r$ form a decreasing filtration of
the space $A^1$ with $\pj=0,\ldots,r         .$ The filtration is preserved by maps
$\pa_V$ and $\pa^*_V,$ hence they induce maps on the associated graded
bundle (we denote them for simplicity of notation by the same symbols as for the full
filtration of $A$). 
We can consider the Kostant Laplacian $\square=gr\,\pa^*_V\,gr\,\pa_V+gr\,\pa_V\,gr\,\pa_V^*.$
The standard Kostant decomposition says that $\Ker\square,$ $\Im gr\,\pa_V^*$ and $\Im gr\,\pa_V$ are complementary subspaces of the
graded bundle $gr\,\cE^i(V)^1.$ In particular, $\square$ is invertible on $\Im gr\,\pa_V^*.$

 Let us consider two consecutive squares with operators
$E_k=d^\Na$ and $E_{k+1}=d^\Na.$ We know that the operator
$G:=\pa^* \circ E_{k+1}\circ E_k =\pa^*_V (R^{d^\Na})$ belongs to $A^1$ and
that the $k$-th  square is commutative iff $G = 0.$ 
If it is not the
case, we shall consider the maximal index $\pj=0$ with the property that 
$G\in A^\pj.$ 

The map $\Phi^{(1)}=-\square^{-1} gr\,(G)$ can be lifted
to a linear algebraic map $\Phi^{(1)}:\cE^k(V)\mapsto\cE^{k+1}(V)$ (e.g., by a choice 
of the Weyl structure) and we shall define the first iteration
$E_k^{(1)}=d^\Na+\Phi^{(1)}.$
Note that the lowest homogeneous component of $E_k^{(1)}$ remains to be $\pa_V$
and that the image of $E_k^{(1)}$ is a subset of $\Im\pa^*.$

Due to 
$$
E_{k+1}\circ E_k^{(1)}-E_{k+1}\circ E_k=
d^\Na \circ \Phi^{(1)} ,
$$
we get 
\begin{eqnarray*}
gr_\pj(\pa^*\circ E_{k+1}\circ E_k^{(1)}))&=&
gr_\pj(G+\pa^*\circ d^\Na \circ\Phi^{(1)})=\\
&=& gr_\pj(G)
-  (gr\,\pa^*)(gr\,\pa_V)(\square^{-1}(gr_\pj(G)))=0.
\end{eqnarray*}
Hence the first order differential operator
$G^{(1)}:=\pa^*\circ E_{k+1}\circ E_k^{(1)}$
belongs to $ A^{\pj+1}.$ 

The same procedure will be repeated inductively.  
If we define
$$
\Phi^{(2)}=-({gr\,\pa^*})\square^{-1}gr_{\pj+1}(G^{(1)})
$$
we can again lift this first order differential operator to a first order differential
operator
$\Phi^{(2)}:\cE^k(V)\mapsto\cE^{k+1}(V)$ and we can define the next iteration by
$$
E_{k}^{(2)}:=E_k^{(1)}+\Phi^{(2)}.
$$
Then
we get 
\begin{eqnarray*}
gr_\pj(\pa_V^*\circ d^\Na\circ E_k^{(2)}))&=&
gr_\pj(G^{(1)}+\pa_V^*\circ d^\Na \circ\Phi^{(2)})=\\
&=& gr_\pj(G^{(1)})
-  (gr\,\pa_V^*)(gr\,\pa_V)({\square}^{-1}(gr_\pj(G^{(1)})))=0.
\end{eqnarray*}
Hence the first order differential operator
$G^{(2)}:=\pa^*_V\circ d^\Na\circ E_k^{(2)}$
belongs to $ A^{\pj+2}.$ 

It is clear that by a finite number of iterations, we shall get the existence part of the 
theorem.

The proof of the uniqueness part of the theorem is similar to the case of Theorem \ref{uniq}.
Suppose that we have two differential operators $\Phi'_k$ and $\Phi''_k$ satisfying the
conditions of the theorem. Their difference $\Phi=\Phi'_k-\Phi''_k$
satisfies $\pa^*_V(d^\Na\circ \Phi)=0.$
 To show that $\Phi=0,$ suppose that $\Phi$ is nontrivial and consider
the biggest $\pj$ such that $(\Phi)^\pj$ is nontrivial.
Then we know that $gr_\pj(d^\Na\circ \Phi)=(\gr\,\pa_V)(gr_\pj(\Phi)),$
hence $(\gr\,\pa_V)(gr_\pj(\Phi))$ is at the same time in  $\Im gr\, \pa_V$
and $\Ker gr\,\pa^*_V,$ hence  $(\gr\,\pa_V)(gr_\pj(\Phi))=0.$
By definition,  $gr_\pj(\Phi)$ belongs also to $\Im \pa^*_V,$ hence
$gr_\pj(\Phi)$ is trivial and we have a contradiction.

As for the second part of the theorem, let us consider two consecutive squares
in the BGG construction induced by $E_k=d^\Na,$ containing operators
$D_k$ and $D_{k+1}.$ If $\Phi_k$ is the deformation constructed above, then
the replacement of $E_k=d^\Na$ by $\tilde{E}_k=d^\Na+\Phi_k$
leads to the same splitting operator $L_k.$ Hence by the first part of the 
theorem,
the $k$-th diagram commutes. Note that the change of the next operator
$E_{k+1}$ will not change the splitting operator $L_{k+1},$
hence the commutativity of the $k$-th diagram is preserved.

Finally, during the construction there were several choices made but due to
the uniqueness of the result, the construction depends only on data of the chosen parabolic geometry.
The same is true
 for the covariant derivative $\nabla.$

\hfill$\square$

\section{Examples.}
 
We want to illustrate in this section general results presented above
by explicit examples showing a form of the prolongation covariant
derivative in some simple situations. Some basic examples in conformal
geometry can be found in \cite{H}. A more comprehensive set of examples
will be treated in \cite{hsss-exmp}.
 
 To calculate the prolongation connection of the first
BGG-operator $ D_0$\ for some tractor bundle $ V=\cG\times_P \bV,$\
we employ the theory of Weyl structures \cite{cap-slovak-weyl},
\cite{cap-slovak-par}.
Both of our examples below will be $|1|$-graded parabolic geometries,
$\g=\g_{-1}\oplus\g_0\oplus\g_1$. Modding out $P_+\cong \g_1$\ of
the parabolic structure bundle $\cG$, one obtains $\cG_0:=\cG/P_+$,
which is a $G_0$-principal bundle over $M$. A splitting $\si:\cG_0\goesto\cG$\
of the canonical projection $\cG\goesto\cG_0$\ is called a \emph{Weyl structure},
and for our geometric structures below this can be identified with
the choice of a \emph{Weyl connection}, 
which is a linear connection $D$\ compatible with the geometry.
Under such a choice, all $P$-associated bundles reduce to $G_0$-associated
bundles, and in particular one gets a decomposition of the tractor bundle $ V$\ which depends on the choice of Weyl structure. The adjoint
tractor bundle $\A M=\cG\times_P\g$\ decomposes into $\A_- M\oplus \A_0 M\oplus \A_1 M$, with $\A_- M\cong TM$\ and $A_+ M\cong T^*M$.
The Lie algebraic action  of $\g$\ on $\bV$\ gives
rise to an action $\bdot$ of $\A M$\ on $ V$, which we can restrict
to $TM$\ and $T^*M$.
The tractor connection $ \na^{\om}$\ can be written as $ \na^{\om}=\del+D+\Rho\bdot$:
the map $\del: V\goesto\Om^1(M, V)$\ is obtained by the action
of $TM\embed \A M$\ on $ V$, and $\Rho\bdot: V\goesto\Om^1(M, V)$\ is
induced by the action of the second slot of the (generalized) Schouten tensor $\Rho\in\ce_{ab}$\ of $D,$ which will be symmetric for our choices of  $D.$ Recall that this decomposition of $ \na^{\om}$\ depends
on the choice of Weyl structure $\si:\cG_0\goesto\cG$\ resp. Weyl connection
$D$.

In our explicit formulas, we  employ abstract index
notation \cite{penrose-rindler-87}: \linebreak $\ce_a=\Om^1(M),\ce^a={\mathfrak{X}}(M)$, multiple
indices are tensor products. Round brackets denote symmetrizations of the
enclosed indices, square brackets denote skew symmetrizations. A subscript zero
takes the trace-free part. 

We are now going to prolong an interesting equation in projective geometry
which has already been treated in \cite{eastwood-matveev} by different
methods, and another equation for Grassmannian structures of type $(2,q),q>2$.
For a more detailed exposition of explicit calculations cf. \cite{mrh-thesis},\cite{H}
and the forthcoming \cite{hsss-exmp}.

\subsection{An example in projective geometry}

Let $M$\ be an orientable manifold of dimension $n$\ endowed with a projective class of linear,
torsion-free connections $[D]$; here $D$\ and $D'$\ are
projectively equivalent if there is a $\Ups\in\ce^1$\ such that $${D'}_a \om_b=D_a\om-\Ups_a\om_b-\Ups_b\om_a,$$
see e.g. \cite{eastwood-notes}.
For simplicity, we will assume that our chosen representatives $D\in[D]$\ preserve a volume form on $TM$.

To define projectively invariant operators we need to employ the
\emph{projective densities}, which are line bundles $\ce[w],w\in\rr$\
associated to the full $\GL(n)$-frame bundle of $TM$\ via
the $1$-dimensional representation 
$$
C\in\GL(n)\mapsto \lvert \det C\rvert^{w\frac{n+1}{n}}\in\rr_+.
$$
We are going to prolong the following projectively invariant operator,
which is written down with respect to a $D\in[D]$, but does not depend on
this choice:
\begin{align}\label{projop}
   D_0&:\ce^{(ab)}[-2]\goesto {\ce_{c\; }^{(ab)}}_0[-2], \\\notag
  \si^{ab}&\mapsto {D_c\si^{ab}}-\frac{1}{n+1}\de_{c}^{(a}D_p\si^{b)p}.
\end{align}
$ D_0$ projects the Levi-Civita derivative of a symmetric two tensor $\si$
to its trace-free part. This operator was discussed in \cite{eastwood-matveev}, where M. Eastwood and V. Matveev showed
that this equation governs the metrizability of a projective class of covariant derivatives.

\subsubsection{The projective structure as a parabolic geometry}
It is a classical result that $(M,[D])$\ is equivalent
to a unique Cartan geometry $(\G,\om)$\
of type $(G,P)=(\SL(n+1),P)$\ with $P$\ the stabilizer of a ray in $\rr^{n+1}$,
see \cite{sharpe},\cite{cap-slovak-par}.

The Lie algebra $\g=\fsl(n+1)$\ is $1$-graded
$\g=\g_{-1}\oplus \g_0\oplus \g_1=\rr^n\oplus \gl(n)\oplus({\rr^n})^*$,
where an element $X\oplus (\al \id+A)\oplus \ph\in\g$\ for $\al\in\rr,A\in\fsl(n)$\ corresponds to the matrix
\begin{align*}
\begin{pmatrix}
  -\al \frac{n}{n+1}  & -\ph \\
  X & \frac{1}{n+1}\al\one_n+A.
\end{pmatrix}
\end{align*}
The actions of $\g_0=\gl(n)\subset\g$\ on
$\g_{-1}=\rr^n$\ and $\g_{1}=({\rr^n})^*$\ are the
standard representation and its dual.

The curvature of the Cartan connection form $\om$\ can be regarded
as an element of $\ce^2(\A M)$, with $\A M=\G\times_P \g$\ the
adjoint tractor bundle, and is written
\begin{align*}
  K=
\begin{pmatrix}
  0 & -A_{ac_1c_2} \\
  0 & {{C_{c_1c_2}}^a}_b
\end{pmatrix}
\end{align*}
with $A$\ the Cotton-York tensor and $C$\ the (projectively invariant)
Weyl curvature (cf. \cite{eastwood-notes}).

$1$-forms and vector fields include into $\A M$\ as
\begin{align*}
\eta_a\in T^*M\mapsto
\begin{pmatrix}
  0  & -\eta_a \\
  0 & 0
\end{pmatrix}
\in \A M,
&
&\xi^a\in TM\mapsto
\begin{pmatrix}
  0  & 0 \\
  \xi & 0
\end{pmatrix}
\in\A M.
\end{align*}

\subsubsection{The operator $ D_0$\ as the first BGG-operator}
Let $ V:=\G\times_P S^2\rr^{n+1}$.
With respect to a choice of Weyl connection $D\in [D]$,
a section $s$ of $ V$\ can be written
\begin{align}\label{vsec}
[s]_D=
\begin{pmatrix}
  \rh \\
  \mu^a \\
  \si^{ab}
\end{pmatrix}
\in
\begin{pmatrix}
   V_2 \\
   V_1 \\
   V_0
\end{pmatrix}
:=
\begin{pmatrix}
  \ce[-2]\\
  \ce^a[-2]\\
  \ce^{(ab)}[-2]
\end{pmatrix}.
\end{align}
We will need that on the first chain spaces the Lie algebra differentials $\del$\ and $\delstar$\ are explicitly given by
\begin{align*}
  &\del\begin{pmatrix}
  \rh \\
  \mu^a \\
  \si^{ab}.
\end{pmatrix}
=
\begin{pmatrix}
  0 \\
  \rh{\de_c}^a \\
  {\de_c}^{(a}\mu^{b)}
\end{pmatrix}
&
&\del\begin{pmatrix}
  \rh_c \\
  {\mu_c}^a \\
  {\si_c}^{ab}.
\end{pmatrix}
=
\begin{pmatrix}
  0 \\
  2{\de_{[c_1}}^a\rh_{c_2]} \\
  2{\de_{[c_1}}^{(a_1}{\mu_{c_2]}}^{a_2)}
\end{pmatrix}.
\\
&\del^*\begin{pmatrix}
  \rh_c \\
  {\mu_c}^a \\
  {\si_c}^{ab}.
\end{pmatrix}
=
\begin{pmatrix}
  -2{\mu_p}^p \\
  -2{\si_p}^{pa} \\
  0
\end{pmatrix}
&
&\del^*\begin{pmatrix}
  \rh_{c_1 c_2} \\
  {\mu_{c_1c_2}}^a \\
  {\si_{c_1c_2}}^{ab}.
\end{pmatrix}
=
\begin{pmatrix}
  2{\mu_{cp}}^p \\
  2{\si_{cp}}^{pa}\\
  0
\end{pmatrix}.
\end{align*}
As bundles with structure group $G_0$, $ V_2, V_1$\ and $T^*M\otimes  V_2$\
are irreducible and are contained
in the image of $\delstar$; $T^*M \otimes  V_1$\ decomposes into
the trace-free part $\im\delstar\cap T^*M\otimes  V_1$\
and the trace part, which lies in the image of $\del$.
The Kostant Laplacian $\lapl$\ acts by
\begin{align*}
\lapl
\begin{pmatrix}
  \rh_{c_1 c_2} \\
  {\mu_{c_1c_2}}^a \\
  {\si_{c_1c_2}}^{ab}
\end{pmatrix}
=
\begin{pmatrix}
  -2n\rh_{c_1 c_2} \\
  -(n+1){\mu_{c_1c_2}}^a \\
  0
\end{pmatrix}
\end{align*}
on $ V$, by multiplication with $-2(n-1)$\ on $T^*M\otimes  V_2$\
and by multiplication with $-n$\ on the trace-free part of
$T^*M\otimes  V_1$.
This is all the algebraic information we need to calculate
the splitting operators and the prolongation.

The tractor connection $ \na^{\om}$ on $ V$ is 
easily calculated with the above actions of $\ce_a$\
and $\ce^a$\ on $ V$\ together with the formula
$ \na^{\om}=\del+D+\Rho\bdot$:
\begin{align*}
   \na^{\om} 
\begin{pmatrix}
  \rh \\
  \mu^a \\
  \si^{ab}
\end{pmatrix}
=
  \begin{pmatrix}
    D_c\rh-2 \Rho_{ca}\mu^a \\
    D_c\mu^a - 2\Rho_{cb}\si^{ab}+\rh{\de_c}^a \\
    D_c\si^{ab}+{\de_c}^{(a}\mu^{b)}
  \end{pmatrix}.
\end{align*}
One calculates that the
first splitting operator
$L_0:\Ga( H_0)\goesto \Ga( V)$\ is given by
\begin{align*}
  \si^{(ab)}\mapsto 
  \begin{pmatrix}
    \frac{1}{n(n+1)}D_pD_q\si^{pq}+\frac{1}{2n}\Rho_{pq}\si^{pq} \\
    -\frac{1}{n+1}D_p\si^{pa} \\
    \si^{ab}
  \end{pmatrix},
\end{align*}
and composition of $ \na^{\om}\circ L_0$\ with projection to the lowest slot
is seen to yield the operator $ D_0$\ of \eqref{projop}.

\subsubsection{Prolongation of $ D_0$}
We calculate the action of the curvature $K\in\Om^2(M,\A M)$:
\begin{align}\label{metricurv}
  K_{c_1c_2}\bdot\begin{pmatrix}0\\0\\ \si^{ab}\end{pmatrix}=
  	\begin{pmatrix}
    -2A_{pc_1c_2}\mu^p \\
    -2A_{pc_1c_2}\si^{pa} 
    + C_{c_1 c_2\;\ p}^{\quad \; a}\mu^p  \\
    2 C_{c_1 c_2\;\ p}^{\;\ \; (a_1}\si^{a_2)p}
    \end{pmatrix}.
\end{align}
Therefore we define
\begin{align*}
  \Phi_1(\begin{pmatrix}0\\0\\\si^{ab}\end{pmatrix}):=
  \begin{pmatrix}
    0 \\
    \bar{\Phi}_1\si \\
    0
  \end{pmatrix}
  :=
-\lapl^{-1}(\delstar(K\bdot\begin{pmatrix}0\\0\\ \si^{ab}\end{pmatrix}))=
  \begin{pmatrix}
  	0 \\
    \frac{2}{n}C_{cp\;\ q}^{\;\ \; a}\si^{pq} \\
    0
  \end{pmatrix}.
\end{align*}
Now the curvature of the modified connection $ \na^{\om}+\Phi_1$\ is
$R=K\bdot+\dcov\Phi_1$
since $(\Phi_1\wedge\Phi_1)(\xi,\eta)$\ vanishes.
For $\xi_1,\xi_2\in {\mathfrak{X}}(M)$\ and $s\in V$\
\begin{align}\label{calc1m}
  &(\dcov\Phi_1)s(\xi_1,\xi_2)=\\&= \notag
\cov_{\xi_1}(\Phi_1(\xi_2)s)-\Phi_1(\xi_2)(\cov_{\xi_1}s)
-\cov_{\xi_2}(\Phi_1(\xi_1)s)+\Phi_1(\xi_1)(\cov_{\xi_2}s)
-\Phi_1([\xi_1,\xi_2])s.
\end{align}
We may expand \eqref{calc1m}\ and write $(\dcov\Phi_1)s$\ as
\begin{align}\label{calc2m}
  \left(\begin{matrix}
   *\\
   \left(\begin{matrix}
 {D}_{\xi_1}\bigl({\bar{\Phi}_1}(\xi_2)\si\bigr)
 -{\bar{\Phi}_1}({\xi_2})\bigl({D}_{\xi_1}\si\bigr)
 -{D}_{\xi_2}\bigl({\bar{\Phi}_1}(\xi_1)\si\bigr)
 +{\bar{\Phi}_1}({\xi_1})\bigl({D}_{\xi_2}\si\bigr)
\\
 -\bar{\Phi}_1([\xi_1,\xi_2])\si
\\
  -{\bar{\Phi}_1}(\xi_2)\del_{\xi_1}\ph
  +{\bar{\Phi}_1}(\xi_1)\del_{\xi_2}\ph
  -{\bar{\Phi}_1}(\xi_2)\del_{\xi_1}\mu
  +{\bar{\Phi}_1}(\xi_1)\del_{\xi_2}\mu        
   \end{matrix}\right)
    \\
    \del_{\xi_1}\bar{\Phi}_1(\xi_2)\si
    -\del_{\xi_2}\bar{\Phi}_1(\xi_1)\si
  \end{matrix}\right),
\end{align}
where we do not  take care about the top component since
it will vanish after an application of $\delstar$.
The lowest component is simply
$\del(\bar{\Phi}_1\si)=-\del\lapl^{-1}\delstar(K\bdot\si).$
Thus $\delstar(Rs)$\ lies in the top slot (i.e., in homogeneity $1$).
So our first adjustment had the effect of moving the
expression $\delstar(Rs)$\ one slot higher.

The new connection $ \na^{\om}+\Phi_1$\ has the following terms
in the middle slot of the curvature $R_{\Phi_1}$:
>From \eqref{calc2m}\ we obtain
the terms $2D_{[c_1}{{{\overline{\Phi}}_1}}_{c_2]}$\ and (via an application
of the algebraic Bianchi identity for $C$), $C_{c_1c_2\; p}^{\; \; \;\ a}\mu^p$.
By \eqref{metricurv}, the contribution of $K\bdot s$\ to
the middle slot is $-2A_{pc_1c_2}\si^{pa} + C_{c_1 c_2\; p}^{\quad \; a}\mu^p$.
In total, we obtain that the action of the curvature $R_{\Phi_1}$\ is
\begin{align*}
  \begin{pmatrix}
    \rh \\
    \mu^a  \\
    \si^{ab}
  \end{pmatrix}
  \mapsto
  \begin{pmatrix}
    * \\
    \frac{2}{n}(D_{[c_1}C_{c_2]p\;\ q}^{\quad \; a})\si^{pq}-2A_{pc_1c_2}\si^{pa} 
    + 2C_{c_1 c_2\;\ p}^{\quad \; a}\mu^p \\
    *
  \end{pmatrix}.
\end{align*}
The entries (*) are irrelevant: the lowest slot is by construction already in the kernel of $\delstar$\ and the highest slot always lies in
$\ker\delstar$.
Now define
\begin{align*}
  \Phi_2(  \begin{pmatrix}
    \rh \\
    \mu^a  \\
    \si^{ab}
  \end{pmatrix}
):=-\lapl^{-1}\delstar(R_{\Phi_1}(\begin{pmatrix}
    \rh \\
    \mu^a  \\
    \si^{ab}
  \end{pmatrix})).
\end{align*}
Using $D_p C_{c_1c_2\;\ a}^{\quad \; p}=(n-2)A_{ac_1c_2}$\ and
trace-freeness of $C,$ we calculate
\begin{align*}
    \Phi_2(  \begin{pmatrix}
    \rh \\
    \mu^a  \\
    \si^{ab}
  \end{pmatrix})=
  \begin{pmatrix}
    -\frac{4}{n} A_{pcq}\si^{pq} \\
    0 \\
    0
  \end{pmatrix}
\end{align*}
and obtain that $\Phi:=\Phi_1+\Phi_2\in\Ga(T^*M\otimes\End( V))$ is 
\begin{align}\label{formulaPsimetric}
  \begin{pmatrix}
    \rh \\
    \mu^a  \\
    \si^{ab}
  \end{pmatrix}
\mapsto
  \frac{2}{n}
  \begin{pmatrix}
    -2 A_{pcq}\si^{pq} \\
    C_{cp\;\ q}^{\;\;\ a}\si^{pq} \\
    0
  \end{pmatrix}.
\end{align}
Now, with $R_{\Phi}$\ the curvature of $\tilde{\na}= \na^{\om}+\Phi$,
one has by construction $\delstar\circ R_{\Phi}=0$.
Thus $\tilde{\na}$ is the prolongation connection for
$({D_c\si^{ab}})_{0} {=}0$.

\subsection{An example in Grassmann geometry}
Let $q\in \bN , q>2$\ and $M$\  be an oriented $2q$-dimensional manifold together with a rank $2$-bundle $\ce_{\al}$\ and a rank $q$-bundle $\ce^{\al'}$.
Assume there is an isomorphism of $TM$\ with $\ce_{\al}\otimes\ce^{\be'}$, which will be fixed. We say that $M$\ together with the identification $TM=\ce_{\al}^{\be'}$ is a \emph{Grassmannian geometry}\ of type $(2,q)$\ if there exists a torsion-free linear connection $D$\ on $TM$\ which is the product of linear connections (again
denoted by $D$) on $\ce_{\al}$\ and $\ce^{\be'}$, see \cite{gover-slovak-quaternionic},
\cite{cap-slovak-par}. The class of all such connections are the Weyl connections of $(M,TM\cong\ce^{\be'}_{\al})$.

We are going to prolong the operator
\begin{align}\label{grassop}
   D_0&:\ce^{[\al'\be']}\goesto ({\ce}_{\ga'}^{\ga[\al'\be']})_0, \\\notag
   D_0(u^{\al'\be'})&=
D^{\ga}_{\ga'}u^{\al'\be'}+\frac{2}{1-q}
\de_{\ga'}^{[\al'}D^{|\ga}_{\tau'}u^{\tau'|\be']}.
\end{align}
Thus, $ D_0(u)$\ is the projection of $Du$\ to its trace-free part.

\subsubsection{Grassmannian structures as parabolic geometries}
Let $G=\SL(n)$, $n=2+q$\ and define $P$\ as the stabilizer of a
two-plane in $(\rr^n)^*$.
Regular, normal and torsion-free parabolic geometries $(\cG,\om)$\ of type $(G,P)$\ are
Grassmannian structures.
In the Cartan-picture, $\ce_{\al}$\ and $\ce^{\al'}$\ are associated
to the $P$-representations $(\rr^p)^*,$\ resp. $(\rr^q)$.

Let $\ST$\ be the standard tractor bundle of $(\cG,\om)$ , i.e., the associated
bundle to the standard representation of $\mbox{\sl SL}(n)$.
Via any Weyl structure $D$, $\ST$\ decomposes into $(\ce^{\al}\oplus\ce^{\al'})$.

The curvature $K\in\ce^2(\mathcal{A} M)=\ce^2(\ST)$ of the Cartan connection
is of the form
\begin{align*}
  K=
  \begin{pmatrix}
    C_{c_1c_2\eta}^{\quad\ \ph} & -A_{pc_1c_2} \\
    0 & {C'}_{c_1c_2\eta'}^{\quad \ph'}
  \end{pmatrix};
\end{align*}
This employs the (generalized) Weyl curvature components
$C\in\Om^2(M,{{\mathfrak{sl}}}(\ce^{\al}))$\ and ${C'}\in\ce^2({\mathfrak{sl}}(\ce^{\al'}))$\ and the
generalized \emph{Cotton-York\ tensor} $A\in\ce^{2}(\ce^1)$ (cf. \cite{gover-slovak-quaternionic}).
Normality of the geometry and torsion-freeness imply
that any possible trace of
$C_{\ga_1'\ga_2'\eta}^{\ga_1\ga_2\ph},{C'}_{\ga_1'\ga_2'\eta'}^{\ga_1\ga_2\ph'}$\
and $A_{\ph'\ga_1'\ga_2'}^{\ph\ga_1\ga_2}$\ vanishes.

\subsubsection{Description of $ D_0$\ as first BGG-operator}
We consider the tractor bundle $ V=\La^2\ST$,
which under choice of a Weyl connection $D$\ decomposes according to
\begin{align*}
  [ V]_{D}=\La^2(\ce^{\al}\oplus\ce_{\al'})=
  \begin{pmatrix}
    \ce^{[\al\be]}\\
    \ce^{\al\be'} \\
    \ce^{[\al'\be']}
  \end{pmatrix}.
\end{align*} 
On the first chain spaces the Lie algebra differentials $\del$\ and $\delstar$\ are
given as follows
(indices within vertical bars are not included in the skew symmetrization):
\begin{align*}
  &\del 
  \begin{pmatrix}
    v^{\al\be}\\
    w^{\al\be'}\\
    u^{\al'\be'}
  \end{pmatrix}
  =
  \begin{pmatrix}
    0 \\
    -\de_{\al'}^{\be'}v^{\al\be}\\
    2\de_{\al'}^{[\be_1'}w^{|\al|\be_2']}
  \end{pmatrix}
&
  &\del 
  \begin{pmatrix}
    v_{\ga'}^{\ga\al\be}\\
    w_{\ga'}^{\ga\al\be'}\\
    u_{\ga'}^{\ga\al'\be'}
  \end{pmatrix}
  =
  \begin{pmatrix}
    0 \\
    2\de_{\ga_2'}^{\be'}v_{\ga_1'}^{\ga_1\ga_2\be}
    -2\de_{\ga_1'}^{\be'}v_{\ga_2'}^{\ga_2\ga_1\be} \\
    2\de_{\ga_1'}^{[\be_1'}w_{\ga_2'}^{|\ga_2\ga_1 | \be_2']}
    -2\de_{\ga_2'}^{[\be_1'}w_{\ga_1'}^{|\ga_1\ga_2 | \be_2']}
  \end{pmatrix}
\\ 
   &\delstar
   \begin{pmatrix}
     v_{\ga'}^{\ga\al\be}\\
     w_{\ga'}^{\ga\al\be'}\\
     u_{\ga'}^{\ga\al'\be'}
   \end{pmatrix}
   =
   \begin{pmatrix}
     -2w_{\tau'}^{[\al_1\al_2]\tau'} \\
     u_{\tau'}^{\al\tau'\be'} \\
     0
   \end{pmatrix}
&
   &\delstar
   \begin{pmatrix}
     v_{\ga_1'\ga_2'}^{\ga_1\ga_2\al\be}\\
     w_{\ga_1'\ga_2'}^{\ga_1\ga_2\al\be'}\\
     u_{\ga_1'\ga_2'}^{\ga_1\ga_2\al'\be'}    
   \end{pmatrix}
   =
  \begin{pmatrix}
    2w_{\ga_1'\tau'}^{\ga_1[\al\be]\tau'}   \\
    -u_{\ga_1'\tau'}^{\ga_1\al\tau'\be'} \\
    0
  \end{pmatrix}.
\end{align*}
The Kostant Laplacian $\lapl=\del\circ\delstar+\delstar\circ\del$\ acts on $[ V]_D$\ via
\begin{align*}
\lapl
  \begin{pmatrix}
    v^{\al\be}\\
    w^{\al\be'}\\
    u^{\al'\be'}
  \end{pmatrix}
=
\begin{pmatrix}
  (2q) v^{\al\be}  \\
  (q-1)w^{\al\be'} \\
  0
\end{pmatrix}.
\end{align*}
The top slot of $\ce^1( V)$\ is
$\ce_c^{[\al\be]}=\ce^{\ga[\al\be]}_{\ga'}$\
and coincides with the image of $\delstar$.
It is irreducible and 
the Kostant Laplacian acts by multiplication with $2(2q-1)$.
The middle slot of $\ce^1( V)$, which is $\ce_c^{\al\be'}$, decomposes
into $\im\del$, which are traces, and the trace-free part
$\im\delstar={\ce_0}_c^{\al\be'}$.
One has that ${\ce_0}^{\ga\al\be'}_{\ga'}=
{\ce_0}^{[\ga \al]\be'}_{\ga'}\oplus{\ce_0}^{(\ga\al)\be'}_{\ga'}$\ and
$\lapl$\ acts by $q$\ on the alternating part
and by $q-2$\ on the symmetric part.

The tractor connection on $ V$\ is
\begin{align*}
   (\na^{\om})^{\ga}_{\ga'}
  \begin{pmatrix}
    v^{\al\be}\\
    w^{\al\be'}\\
    u^{\al'\be'}
  \end{pmatrix}
  =
  \begin{pmatrix}
    D^{\ga}_{\ga'}v^{\al\be}+2\Rho_{\ga'\tau'}^{\ga[\al}w^{\be]\tau'} \\
    D^{\ga}_{\ga'}w^{\al\be'}-\de^{\al}_{\ga'}v^{\ga\be'}+\Rho^{\ga\al}_{\ga'\ta'}u^{\be'\tau'}\\
    D^{\ga}_{\ga'}u^{\al'\be'}+2\de^{[\al'}_{\ga'}w^{|\ga|\be']}
  \end{pmatrix}.
\end{align*}
The first BGG-splitting operator $L_0:\ce^{(\al'\be')}\goesto\Ga( V)$\ is computed
\begin{align*}
  L_0(u^{\al'\be'})=
  \begin{pmatrix}
    \frac{1}{2q}\Rho_{\tau_1'\tau_2'}^{\al\be}u^{\tau_1'\tau_2'}-\frac{1}{1-q}D_{\tau_1'}^{[\al}D_{\tau_2'}^{\be]}u^{\tau_1'\tau_2'} \\
    \frac{1}{1-q}D_{\tau'}^{\al}u^{\tau'\be'} \\
    u^{\al'\be'}
  \end{pmatrix},
\end{align*}
and composition of $ \na^{\om}\circ L_0$\ with projection to the lowest slot
is seen to yield our operator \eqref{grassop}.

\subsubsection{Prolongation of $ D_0$}
For a section $s$\ of $ V$\ one first computes
$K\bdot s\in\ce^2( V)$,  which
is then mapped by $\delstar$\ into $\ce^1( V)$,
\begin{align}\label{delstar-K-grass}
  \delstar(K\bdot
  \begin{pmatrix}
    v^{\al\be}\\
    w^{\al\be'}\\
    u^{\al'\be'}
  \end{pmatrix}
)
  &=
  \begin{pmatrix}
    2C_{\ga_1'\ph'\eta}^{\ga_1[\al\be]}w^{\eta\ph'}
    +2A_{\eta'\ph'}^{[\al|\ga_1|\be]}u^{\eta'\ph'}
    \\
    -2{C'}_{\ga_1'\ph'\eta'}^{\ga_1\al\be'}u^{\ph'\eta'}
    \\
    0
  \end{pmatrix}.
\end{align}
The first deformation map $\Phi_1$\ is defined by
$\Phi_1=-\lapl^{-1}\circ\delstar\circ K\bdot$,
\begin{align*}
  \Phi_1(\begin{pmatrix} 0 \\ 0 \\ u^{\al'\be'}\end{pmatrix})=
  \begin{pmatrix}
    0\\
  \frac{2}{q}
  {C'}_{\ga_1'\ph'\eta'}^{[\ga_1\al]\be'}u^{\ph'\eta'}
  +\frac{2}{q-2}
  {C'}_{\ga_1'ph'\eta'}^{(\ga_1\al)\be'}u^{\ph'\eta'}
  \\
  0
\end{pmatrix}.
\end{align*}
Now we need to calculate $\delstar$\ of the change in curvature resulting from $\Phi_1$, which is just $\delstar\circ\dcov\Phi_1$,
since one quickly sees that $\delstar \circ {\Phi_1}_{[c_1}\del_{c_2]}=0.$
Both indices of a section $w^{\al\be'}$\ are contracted into $C$
and the trace taken by $\delstar$\ vanishes by trace-freeness of $C$, $C'$.
Therefore we are only interested in the differential components of $\dcov \Phi_1$ given by
\begin{align*}
\begin{pmatrix}
0
\\
\begin{pmatrix}
2(\frac{1}{q}+\frac{1}{q-2})
 D_{\ga_1'}^{\ga_1}{C'}_{\ga_2'\ph'\eta'}^{\ga_2\al\be'}u^{\ph'\eta'}
 -2(\frac{1}{q}-\frac{1}{q-2})
 D_{\ga_1'}^{\ga_1}{C'}_{\ga_2'\ph'\eta'}^{\al\ga_2\be'}u^{\ph'\eta'}
 \\
 -2(\frac{1}{q}+\frac{1}{q-2})
 D_{\ga_2'}^{\ga_2}{C'}_{\ga_1' \ph' \eta'}^{\ga_1\al\be'}u^{\ph'\eta'}
 +2(\frac{1}{q}-\frac{1}{q-2})
 D_{\ga_2'}^{\ga_2}C_{\ga_1'\ph'\eta'}^{\al\ga_1\be'}u^{\ph'\eta'}
\end{pmatrix}
\\
0
\end{pmatrix}.
\end{align*}
Applying $\delstar$\ we obtain the top slot contribution
\begin{align}
  \label{delstar-na-grass}
     -4(\frac{1}{q}+\frac{1}{q-2})
 D_{\tau'}^{[\al}{C'}_{\ga_1' \ph' \eta'}^{|\ga_1|\be]\tau'}u^{\ph'\eta'}
 +4(\frac{1}{q}-\frac{1}{q-2})
 D_{\tau'}^{[\al}C_{\ga_1'\ph'\eta'}^{\be]\ga_1\tau'}u^{\ph'\eta'}
\end{align}
Adding the contributions of the top slot of \eqref{delstar-K-grass}\ and \eqref{delstar-na-grass} 
(after multiplication by $-\frac{1}{2(2q-1)}$) to the modification  map $\Phi_1$,
we obtain the full modification map
\begin{align*}
  \Phi
  \begin{pmatrix}
    v^{\al\be}\\
    w^{\al\be'}\\
    u^{\al'\be'}       
  \end{pmatrix}
  =
  \begin{pmatrix}
    \frac{1}{2q-1}
    \begin{pmatrix}
   2(\frac{1}{q}+\frac{1}{q-2})
 D_{\tau'}^{[\al}{C'}_{\ga_1' \ph' \eta'}^{|\ga_1|\be]\tau'}u^{\ph'\eta'}
 -2(\frac{1}{q}-\frac{1}{q-2}) D_{\tau'}^{[\al}C_{\ga_1'\ph'\eta'}^{\be]\ga_1\tau'}u^{\ph'\eta'}
\\
 -C_{\ga_1'\ph'\eta}^{\ga_1[\al\be]}w^{\eta\ph'}
  +A_{\eta'\ph'}^{[\al|\ga_1|\be]}u^{\eta'\ph'}
    \end{pmatrix}
\\
\begin{pmatrix}
  \frac{2}{q}
  {C'}_{\ga_1'\ph'\eta'}^{[\ga_1\al]\be'}u^{\ph'\eta'}
  +\frac{2}{q-2}
  {C'}_{\ga_1'\ph'\eta'}^{(\ga_1\al)\be'}u^{\ph'\eta'}  
\end{pmatrix}
\\
0
  \end{pmatrix}.
\end{align*}
$\tilde{\na}= \na^{\om}+\Phi$ is then the prolongation connection of the system
$(D^{\ga}_{\ga'}u^{\al'\be'})_0 {=}0$.

\subsection{The case of infinitesimal automorphisms.}
Let $\A M$\ be the adjoint tractor bundle of a regular parabolic
geometry $(\cG,\om)$\ over $M$\ and $ \na^{\om}$\ the adjoint tractor connection.
In \cite{Cap}\ it was shown that parallel sections of the
connection 
\begin{align}\label{solcap}
  \tilde\na s= \na^{\om} s+\ka(\Pi(s),\cdot)
\end{align}
are in 1-1-correspondence with infinitesimal automorphisms of $(\cG,\om),$
where $\Pi$\ is the natural projection $\Pi:\A M\goesto TM$.
It shows that it is interesting to consider the first BGG-operator $\tilde D_0$\ obtained
by $\tilde\na$. If the parabolic geometry $(\cG,\om)$\
is normal, the curvature of $\tilde\na$\ lies in the kernel of $\delstar_{\A M}$.
Therefore, exactly as in corollary \ref{cor-prolong}, one sees that
$\Pi_0:\A M\goesto H_0$\ and $\tilde L_0:H_0\goesto \A M$\ are inverse
isomorphisms between the space of parallel sections of $\tilde\na$\ and
the kernel of $\tilde D_0.$ Thus, the operator $\tilde D_0$\ describes
the infinitesimal automorphisms of $(\cG,\om)$\ and is
automatically prolonged by $\tilde \na$.

It is shown that if the parabolic geometry is also torsion-free or $1$-graded,
one has that $\delstar_{\A M}\ka=0$; i.e., for every $s\in\A M$\ one has
$\delstar_{\A M} \ka(\Pi(s),\cdot)=0$. But in the torsion-free case,
the map $\xi\mapsto \ka(\Pi(s),\xi)$\ is
evidently homogeneous of degree $\geq 1$. Therefore, if
we know that $H_1(\g_-,\g)$\ sits in homogeneity $\leq 0$, we see
that $\xi\mapsto \ka(\Pi(s),\xi)$\ lies in $\im\delstar_{\A M}$.

Thus we have:
\begin{thm}
	Let $(\cG\goesto M,\om)$\ be a torsion-free, normal parabolic geometry
	with $H_1(\g_-,\g)$\ concentrated in homogeneity $\leq 0$.
	Then $\tilde\na$\ from \eqref{solcap}\ coincides with the normalized
	tractor connection on $\A M$.
In particular, the usual first BGG-operator $D_0$\
coincides with $\tilde D_0$\ and thus describes infinitesimal automorphisms.
\end{thm}

We note that the homogeneity condition on $H_1(\g_-,\g)$\ is satisfied
for all parabolic geometries of type $(G,P)$ with $\g$\ simple and
$(G,P)$\ not corresponding to projective structures or contact projective structures.

\section{Acknowledgement.}
The results presented in the paper are inspired very much by construction achieved by A. \v Cap
in a description of infinitesimal automorphisms for parabolic geometries.
The authors would like to thank him for a lot of helpful discussions.
The first author was supported by the IK I008-N funded by the
University of Vienna and project P 19500-N13 of the "Fonds zur
F\"{o}rderung der wissenschaftlichen Forschung" (FWF).
The second and the third authors were supported by the 
institutional grant MSM 0021620839 and
by the grant GA CR 201/08/397.
 The last author was
supported from Basic Research Center no.\ LC505 (Eduard \v{C}ech
Center for Algebra and Geometry) of the Ministry of Education of Czech
Republic.


\begin{thebibliography}{XX}

\bibitem{BCEG}  Branson T.,  \v Cap A., Eastwood M. G., Gover A. R.: Prolongations of geometric overdetermined systems, Int. J. Math.\ {\bf 17} (2006), 641?664.

\bibitem{Cap}  \v Cap A.: Infinitesimal Automorphisms and Deformations of Parabolic Geometries,
J. Eur. Math. Soc.\ {\bf 10}(2008) 415-437.

\bibitem{cap-slovak-weyl} \v Cap A., Slov\'ak J.: Weyl structures for parabolic geometries, Math. Scand.\ {\bf  93} (2003), 53-90.

\bibitem{CD} Calderbank D., Diemer T., Differential invariants and
  curved Bernstein-Gelfand-Gelfand sequences,   J.\ Reine Angew.\ Math.\
  {\bf 537} (2001), 67-103.

\bibitem{cap-slovak-par} \v Cap A., Slov\'ak J.: Parabolic Geometries I: Background and General Theory, AMS, Mathematical Surveys and Monographs, 2009, 628 pp.

\bibitem{CSS} \v{C}ap A., Slov\'{a}k J. and Sou\v{c}ek V.,
  {Bernstein-Gelfand-Gelfand sequences}, Ann.\ Math.\ {\bf 154} (2001),
  97--113.

\bibitem{DT} Dunajski M., Tod P.: Four dimensional metrics conformal to K\"ahler, preprint, 
arXiv:0901.2261, to appear in Mathematical Proceedings of the Cambridge Philosophical Society. 

\bibitem{eastwood-notes} Eastwood M. G.: Notes on Conformal Differential Geometry, Supp. Rend. Circ. Matem. Palermo {\bf 43} (1996), 57-76. 

\bibitem{eg} Eastwood M. G., Gover A. R.: Prolongation on contact manifolds, preprint, arXiv:0910.5519.

\bibitem{eastwood-matveev} Eastwood M. G., Matveev V.: Metric Connections in Projective Differential Geometry, Symmetries and Overdetermined Systems of Partial Differential Equations,
IMA Vol. Math. Appl., Vol. 144, (2007), pp. 339-351. 

\bibitem{E} Eastwood M. G.: Notes on projective differential geometry, in Symmetries and overdetermined systems of aprtial differential equations, IMA Vol. Math. Appl., {\bf 144}, Springer, New York, 2008, 41-60.


\bibitem{GoEinst} Gover A. R.: Laplacian operators and 
  $Q$-curvature on conformally Einstein manifolds, 
  Math.\ Ann.\ {\bf  336}  (2006),   311-334.

\bibitem{GSS} Gover A. R., Somberg P., Sou\v{c}ek V.: 
  Yang-Mills detour complexes and conformal geometry, 
  Commun.\ Math.\ Phys.\ {\bf 278} (2008), 307-327. 

\bibitem{GS} Gover A. R., \v{S}ilhan J.:
  The conformal Killing equation on forms - prolongations 
  and applications, Diff. Geom. and its Appl. {\bf  26} (2008), 244-266. 
  
\bibitem{gover-slovak-quaternionic} Gover A. R., Slov\'ak J.: Invariant Local Twistor Calculus for Quaternionic Structures and Related Geometries, Journal of Geometry and Physics {\bf 32} (1999), 14-56. 


\bibitem{H} Hammerl M.: Invariant prolongation of BGG-operators in conformal geometry,  Arch.\ Math.\ {\bf  44} (2008), 367-384.  

\bibitem{mrh-thesis} Hammerl M.: Natural prolognations of BGG operators, Thesis, University of Vienna, 2009.

\bibitem{hsss-exmp} Hammerl M., Somberg P., Sou\v cek V., \v{S}ilhan J.: On a new normalization of tractor covariant derivatives: Examples, preprint.

\bibitem{penrose-rindler-87} Penrose R., Rindler W.: Spinors and space-times, Vol. 1,2,
Cambridge University Press, 1982, 1984.

\bibitem{sharpe} Sharpe R. W.: Differential Geometry, Graduate Texts in Mathematics 166, Springer?Verlag, 1997.
\end{thebibliography}
\end{document}